\documentclass[a4paper, 12pt]{article}
\usepackage{amsmath} \usepackage{amsfonts} \usepackage{amssymb}\usepackage{latexsym}
\usepackage[english]{babel}
\usepackage{amscd}
\usepackage[dvips,final]{graphics}
\usepackage{locengli,math}
\usepackage{graphicx, Addaletree}
\usepackage[all]{xy}
\usepackage{latexsym}
\usepackage{eufrak}

\begin{document}
\begin{center}
\textbf{\LARGE{\textsf{Hochschild two-cocycles and the good triple $(As,Hoch,Mag^\infty)$
}}}~\footnote{
\textsf{Email}: ph$\_$ler$\_$math@yahoo.com}
\vskip2cm
\large{Philippe {\sc Leroux}}
\end{center}
\vskip2cm
\noindent
{\bf Abstract:}
Hochschild two-cocycles play an important role in the deformation \`a la Gerstenhaber of associative algebras.
The aim of this paper is to introduce the category of Hoch-algebras whose objects are associative algebras equipped with an extra magmatic operation $\succ$ verifying the Hochschild two-cocycle relation:
$$ \mathcal{R}_2: \ \ (x \succ y)*z+ (x*y)\succ z= x\succ (y*z)+ x*(y\succ z).$$
The free Hoch-algebra over a $K$-vector space is given in terms of planar rooted trees and the triples of operads
$(As,Hoch, Mag^\infty)$ endowed with the infinitesimal relations are shown to be good. We then obtain an equivalence of categories between connected infinitesimal $Hoch$-bialgebras and $Mag^\infty$-algebras.

\noindent
\textbf{Notation}:
In the sequel $K$ is a  field. We adopt Sweedler notation for the binary cooperation $\Delta$ on a $K$-vector space $V$ and set $\Delta(x)= x_{(1)} \otimes x_{(2)}$. For a $K$-vector space $V$, we set $\bar{T}(V):=\bigoplus_{n>0} V^{\otimes n}$.

\section{Introduction}
The well-known Poincar\'e-Birkhoff-Witt and the Cartier-Milnor-Moore theorems together
can be rephrased as follows:
\begin{theo}(CMM-PBW)
For any cocommutative (associative) bialgebra $\mathcal{H}$, $Com^c-As$-bialgebra for short, the following is equivalent.
\begin{enumerate}
 \item{$\mathcal{H}$ is connected;}
\item{$\mathcal{H}$ is isomorphic to $U(Prim \ \mathcal{H})$ as a bialgebra;}
\item{$\mathcal{H}$ is isomorphic to $Com^c(Prim \ \mathcal{H})$ as a coalgebra,}
\end{enumerate}
where $U$ is the usual enveloping functor and $ Prim \ \mathcal{H}$ the usual Lie algebra of the primitive elements of $\mathcal{H}$.
\end{theo}
\noindent
In the theory developed by J.-L. Loday \cite{GB}, this result is rephrased by saying that the triple of operads $(Com, As, Lie)$, endowed with the usual Hopf relation, is good, where $Com$, $As$, and $Lie$ stand respectively for the operads of commutative, associative and Lie algebras.
Other good triples of operads equipped with other relations than the usual Hopf one, have been found since. A summary can be found in \cite{GB}, see also \cite{Ltrip, infini} for other examples.

It has been shown in \cite{infini} that the triple of operads $(As, Dipt, Mag^\infty)$ endowed with the semi-infinitesimal relations is good. The operad $Dipt$ is related to dipterous algebras which are associative algebras equipped with an extra left module on themselves, see also \cite{LodRon}, and $Mag^\infty$ is related to $Mag^\infty$-algebras, i.e., $K$-vector spaces having one $n$-ary (magmatic) generating operation for each integer $n>1$.
We then obtained that the category of
connected infinitesimal dipterous bialgebras, $As^c-Dipt$-bialgebras for short, was equivalent to the category of $Mag^\infty$-algebras.
In this paper, we propose another equivalence of category involving $Mag^\infty$: the category of
connected infinitesimal $Hoch$-bialgebras is equivalent to the category of $Mag^\infty$-algebras.

In Section~2, we introduce $Hoch$-algebras and give an explicit construction of the free
$Hoch$-algebra over a $K$-vector space. In Section~3, we introduce the notion of (connected) infinitesimal $Hoch$-bialgebras. In Section~4 we prove the announced equivalence of categories. In Section~5, we deal with unital $Hoch$-algebras and close by Section~6 with two other good triples involving the operad $Hoch$.
\section{The free $Hoch$-algebra}
A $Hoch$-algebra $G$ is a $K$-vector space equipped with an associative operation $*$ and a magmatic operation $\succ$ verifying:
$$\mathcal{R}_2: (x \succ y)*z+ (x*y)\succ z= x\succ (y*z)+ x*(y\succ z),$$
for all $x,y,z\in G$.
Let $V$ be a $K$-vector space. The free $Hoch$-algebra over $V$ is defined as follows.
It is equipped with a linear map $i: V \rightarrow Hoch(V)$ and
for any $Hoch$-algebra $G$ and any linear map $f: V \rightarrow G$, there exists a unique $Hoch$-algebra morphim $\phi: Hoch(V)\rightarrow G$ such that $\phi \circ i=f$. We now give an explicit construction of the free $Hoch$-algebra over a $K$-vector space.

\noindent
Denote by $T_n$ the set of rooted planar trees (degrees at least 2) with $n$ leaves. The cardinalities of $T_n$ are registered under the name \textit{A001003 little Schroeder numbers} of the Online Encyclopedy of Integer Sequences. For $n=1,2,3$, we get:
$$ T_1=\{\ \vert \ \}, \ T_2=\{\ \treeA \ \}, \ T_3=\{\ \treeAB, \treeBA, \ \treeM \ \}.$$
Define grafting operations
by:
$$ [\cdot, \ldots,\cdot]:T_{n_1} \times \ldots \times T_{n_p} \rightarrow T_{n_1 + \ldots + n_p}, \ \ (t_1, \ldots, t_p)\mapsto [t_1, \ldots, t_p]:=t_1\vee \ldots \vee t_p,$$
where the tree $t_1\vee \ldots \vee t_p$ is the tree whose roots of the $t_i$ have been glued together and a new root has been added. Observe that any rooted planar tree $t$ can be decomposed in a unique way via the grafting operation as
$t_{1} \vee \ldots \vee t_{p}$.
Set $T_{\infty}:= \bigoplus_{n>0} \ KT_n$.
Define over $\bar {T}(T_{\infty})$, the following binary operations, first on trees, then by bilinearity:
$$ Concatenation: \ \ (t_1 \ldots t_p) * (s_1 \ldots s_q):= t_1 \ldots t_ps_1 \ldots s_q,$$
$$ (t_1 \ldots t_p) \succ (s_1 \ldots s_q):= \sum_{k=1}^q \sum_{i=0}^{p-1} \ t_1 \ldots t_{p-(i+1)}[t_{p-i},\ldots,t_p, s_1, \ldots, s_k]s_{k+1}\ldots s_q.$$
For instance we get:
$$ | \ | \ |\succ |:= | \ | \ \treeA + | \ \treeM +\treeCor $$
$$ | \succ | \ \treeA:= \treeA  \ \treeA+ \treeMA$$
\begin{theo}
The $K$-vector space $\bar {T}(T_{\infty})$ endowed with the operations $*$ and $\succ$ is the free $Hoch$-algebra over $K$.
\end{theo}
\Proof
Let $x:=x_1 \ldots x_m$, $y:=y_1 \ldots y_n$ and $z:=z_1\ldots z_p$. We get:
\begin{eqnarray*}
(x\succ y)*z + (x*y)\succ z &=& \sum_{k=1}^n \sum_{i=0}^{m-1} \ x_1 \ldots x_{m-(i+1)}[x_{m-i},\ldots,x_m, y_1, \ldots, y_k]y_{k+1}\ldots y_nz_1\ldots z_p\\
& + &  \sum_{k=1}^p \sum_{i=0}^{n-1} \ x_1 \ldots x_m y_1 \ldots y_{n-(i+1)}[y_{n-i},\ldots,y_n, z_1, \ldots, z_k]z_{k+1}\ldots z_p\\
&+ & \sum_{k=1}^p \sum_{i=0}^{m-1} \ x_1 \ldots x_{m-(i+1)}[x_{m-i},\ldots,x_m, y_1, \ldots, y_n, z_1 \ldots, z_k]z_{k+1}\ldots z_p
\end{eqnarray*}

\begin{eqnarray*}
x\succ (y*z) + x*(y \succ z) &=&
\sum_{k=1}^n \sum_{i=0}^{m-1} \ x_1 \ldots x_{m-(i+1)}[x_{m-i},\ldots,x_m, y_1, \ldots, y_k]y_{k+1}\ldots y_nz_1 \ldots z_p \\
&+& \sum_{k=1}^p \sum_{i=0}^{m-1} \ x_1 \ldots x_{m-(i+1)}[x_{m-i},\ldots,x_m, y_1, \ldots, y_n,z_1, \ldots, z_k]z_{k+1}\ldots z_p\\
&+& \sum_{k=1}^p \sum_{i=0}^{n-1} \ x_1 \ldots x_m y_1 \ldots y_{n-(i+1)}[y_{n-i},\ldots,y_n, z_1, \ldots, z_k]z_{k+1}\ldots z_p,
\end{eqnarray*}
showing that
$$ (x \succ y)*z+ (x*y)\succ z= x\succ (y*z)+ x*(y\succ z),$$
holds for all forests of planar rooted trees $x,y,z$.
Observe that any rooted planar tree $t:= [t_1, \ldots, t_n]$
can be rewritten as:
$$ t= (t_1*(t_2 \ldots t_{n-1}))\succ t_n - t_1*((t_2 \ldots t_{n-1})\succ t_n)$$
Let $G$ be a $Hoch$-algebra and $g\in G$ and $f: K \rightarrow G$ be a linear map. Consider the embedding $i: K \hookrightarrow \bar {T}(T_{\infty})$ defined by $i(1_K):=|$ and define by induction the map $\phi: \bar {T}(T_{\infty}) \rightarrow G$ as follows:
$$ \phi(|)=g, $$
$$\phi(t_1 \ldots t_n)=\phi(t_1)*_G \phi(t_2)*_G \ldots *_G \phi(t_{n}),$$
$$\phi(t)= (\phi(t_1)*_G(\phi(t_2) \ldots \phi(t_{n-1})))\succ_G \phi(t_n) - \phi(t_1) *_G ((\phi(t_2) \ldots \phi(t_{n-1}))\succ_G \phi(t_n)) $$
for any $t:= [t_1, \ldots, t_n]$
and extend $\phi$ by linearity.
By construction, $\phi$ is a morphism of associative algebras. Using
the fact that:
$$  (x*y)\succ z- x*(y\succ z)= x\succ (y*z)-(x \succ y)*z,$$
and changes of indices in the involving sums, one shows that $\phi$ is also a morphism for the magmatic operations.
It is then the only $Hoch$-morphism such that $\phi \circ i= f$.
\eproof

\noindent
As the operad $Hoch$ is regular, the following holds.
Let $V$ be a $K$-vector space. The free $Hoch$-algebra over $V$ is the $K$-vector space:
$$ Hoch(V):= \bigoplus_{n>0} {Hoch}_n \otimes V^{\otimes n}, $$
with  $Hoch(K):= \bigoplus_{n>0} {Hoch}_n\simeq \bar {T}(T_{\infty})$ (hence ${Hoch}_n$ is explicitely described in terms of forests of rooted planar trees) equipped with the operations $*$ and $\succ$ defined as follows:
$$ ((t_1\ldots t_n) \otimes \omega) * ((s_1\ldots s_p) \otimes \omega')= (t_1\ldots t_ns_1\ldots s_p) \otimes \omega \omega', $$
$$ ((t_1\ldots t_n) \otimes \omega) \succ ((s_1\ldots s_p) \otimes \omega')= ((t_1\ldots t_n)\succ( s_1\ldots s_p)) \otimes \omega \omega',$$
for any $\omega \in V^{\otimes n}, \omega'\in V^{\otimes p}$.
The embedding map $i: V \hookrightarrow Hoch(V)$ is defined by: $v \mapsto \vert \otimes v$.

\noindent
Since the generating function associated with the Schur functor $\bar {T}$ is $f_{\bar {T}}(x):=\frac{x}{1-x}$ and with the Schur functor $T_{\infty}$ is
$f_{T_{\infty}}(x):=\frac{1 + x - \sqrt(1 - 6x + x^2)}{4}=x+x^2+3x^3+11x^4+45x^5 + \ldots$,
the generating function of the operad $Hoch$ is $f_{\bar {T}}\circ f_{T_{\infty}}$, that is:
$$ f_{Hoch}(x):=\frac{1+x-\sqrt{1-6x+x^2}}{3-x+\sqrt{1-6x+x^2}}=x+2x^2+6x^3+22x^4+\ldots.$$
The sequence $(1,2,6,22,90,\ldots)$ is registered as $A006318 $ under the name \textit{Large Schroeder numbers} on the Online Encyclopedy of Integer Sequences.

\section{Infinitesimal $Hoch$-bialgebras}

By definition, an infinitesimal $Hoch$-bialgebra (or an $As^c-Hoch$-bialgebra for short) $(\mathcal{H},* \succ,\Delta)$ is a $Hoch$-algebra equipped with a coassociative coproduct $\Delta$ verifying the following so-called nonunital infinitesimal relations:
$$ \Delta(x \succ y):= x_{(1)}\otimes (x_{(2)} \succ y) + (x \succ y_{(1)}) \otimes y_{(2)} + x\otimes y.$$
$$ \Delta(x * y):= x_{(1)}\otimes (x_{(2)} * y) + (x * y_{(1)}) \otimes y_{(2)} + x\otimes y.$$
It is said to be connected when $\mathcal{H}=\bigcup_{r>0} F_r\mathcal{H}$ with
the filtration $(F_r\mathcal{H})_{r>0}$ defined as follows:
$$ \textrm{(The primitive elements)} \ \ F_1\mathcal{H}:=Prim \ \mathcal{H}=\ker \Delta,$$
Set $\Delta^{(1)}:= \Delta$ and
$\Delta^{(n)}:= (\Delta \otimes id_{n-1})\Delta^{(n-1)}$ with
$id_{n-1}= \underbrace{id \otimes \ldots \otimes id}_{times \ n-1}$. Then,
$$ F_r\mathcal{H}:=\ker \ \Delta^{(r)}.$$

\begin{theo}
\label{h1}
Let $V$ be a $K$-vector space. Define on $ Hoch(V)$,
the free $Hoch$-algebra over $V$,
the cooperation $\Delta: Hoch(V) \rightarrow Hoch(V) \otimes Hoch(V)$ recursively as follows:
$$ \Delta(i(v)):=0, \ \textrm{for all} \ v \in V,$$
$$ \Delta(x \succ y):= x_{(1)}\otimes (x_{(2)} \succ y) + (x \succ y_{(1)}) \otimes y_{(2)} + x\otimes y.$$
$$ \Delta(x \star y):= x_{(1)}\otimes (x_{(2)} \star y) + (x \star y_{(1)}) \otimes y_{(2)} + x\otimes y,$$
for all $x,y \in Hoch(V)$.
Then $(Hoch(V), \Delta)$
is a connected infinitesimal $Hoch$-bialgebra.
\end{theo}
\Proof
This result can be proved by hand or can be seen as a corollary of the Theorem~\ref{good} in the next section.
\eproof

\section{A good triple of operads}
It can be usefull to have the following result when searching for good triples.
\begin{lemm}
\label{identif}
Let $\mathcal{C},\mathcal{A},\mathcal{Z},\mathcal{Q}$ and $Prim$ be operads. Suppose the triples
of operads $(\mathcal{C},\mathcal{A}, \ Prim)$ and $(\mathcal{C},\mathcal{Z}, Vect)$  equipped with the same compatibility relations to be good. Suppose $\mathcal{A}=\mathcal{Z} \circ \mathcal{Q}$ then $Prim=\mathcal{Q}$.
\end{lemm}
\Proof
Since $(\mathcal{C},\mathcal{Z}, Vect)$ is good, the notion of $\mathcal{C}^c-\mathcal{Z}$-bialgebra has a meaning and the following is equivalent:
\begin{enumerate}
\item{The $\mathcal{C}^c-\mathcal{Z}$-bialgebra $\mathcal{H}$ is connected.}
\item{As $\mathcal{Z}$-algebras, $\mathcal{H}$ is isomorphic to the free $\mathcal{Z}$-algebra over its primitive elements.}
\item{As $\mathcal{C}^c$-coalgebras, $\mathcal{H}$ is isomorphic to the cofree $\mathcal{C}^c$-coalgebra over its primitive elements.}
\end{enumerate}
As $(\mathcal{C},\mathcal{A}, \ Prim)$ is good, the isomorphism of Schur functors $\mathcal{A}\simeq \mathcal{C}^c \circ Prim$ holds. Therefore, if $V$ is a $K$-vector space, then $\mathcal{A}(V)= \mathcal{C}^c (Prim (V))$ and by hypothesis $\mathcal{A}(V)=\mathcal{Z}(\mathcal{Q}(V))$. Hence $\mathcal{A}(V)$ is a connected $\mathcal{C}^c-\mathcal{Z}$-bialgebra. Consequently, $Prim = \mathcal{Q}$.
\eproof

\noindent
\begin{theo}
\label{good}
The triple of operads $(As, Hoch, Mag^{\infty})$ endowed with the infinitesimal relation is good.
\end{theo}
\Proof
Fix an integer $n>0$.
By $[n]-Mag$ we mean the regular binary operad generated by $n$ magmatic (binary) operations.
In \cite{Ltrip} Thm. 4.4 (and 4.5), it has been shown that for each integer $n>0$ the triples of operads
$(As, [n]-Mag, Prim \ [n]-Mag)$ endowed with the infinitesimal relations were
good. For $n=2$,
the operadic ideal $J$ generated by the primitive operations:
$$*(\succ \otimes  id )+ \succ (* \otimes  id )- \succ(id \otimes  *) - *(id \otimes \succ ),$$
$$*(* \otimes  id )- *(id \otimes * ),$$
yields another good triple of operads $(As, [2]-Mag/J ,  Prim \ ([2]-Mag/J))$ (cf. Prop 3.1.1 \cite{GB} on quotient triples), which turns out to be
the triple $(As, Hoch, Prim \ Hoch)$.
As $(As,As,Vect)$ endowed with the infinitesimal relation is good
(cf \cite{LodRon}) and
$Hoch=As \circ Mag^\infty$ using Section~2, we get $Prim \ Hoch=Mag^\infty$ by using Lemma~\ref{identif}.
\eproof

\noindent
We then obtain another equivalence of categories involving the operad $Mag^\infty$.
\begin{coro}
The category of connected infinitesimal $As^c-Hoch$-bialgebras and the category of $Mag^\infty$-algebras are equivalent.
$$ \{\textrm{conn.} \ As^c-Hoch-bialg.\} \underset{Primitive}{\overset{U}{\leftrightarrows}} \{Mag^\infty-alg.\},$$
where $U$ and $Primitive$ are respectively the universal enveloping functor and the primitive functor.
\end{coro}
\Proof
Apply Thm. 2.6.3. \cite{GB}.
\eproof

\NB
The functor $Primitive$ is obviously given as follows. If $(\mathcal{H},*, \succ)$ is a connected infinitesimal $Hoch$-bialgebra, then for all integer $n>1$ and for all primitive elements $x_1, \ldots,x_n \in \mathcal{H}$, the element:
$$ [x_1, \ldots,x_n]_n:=(x_1* \ldots * x_{n-1})\succ x_n-x_1*((x_2* \ldots * x_{n-1})\succ x_n),$$
will be primitive.
The functor $U$ acts as follows. Let $(M, ([,\ldots,]_n)_{n>1})$ be a $Mag^\infty$-algebra with the $[,\ldots,]_n$ being its generating $n$-ary operations. Then $U(M)$ is given by $Hoch(M)/ \sim$, where the equivalence relation $\sim$ consists in identifying,
$$(x_1* \ldots * x_{n-1})\succ x_n-x_1*((x_2* \ldots * x_{n-1})\succ x_n),$$
with $[x_1, \ldots,x_n]_n$,
for all $x_1 \ldots,x_n \in M$.

\section{Extension to a unit}
Unital $Hoch$-algebras are $Hoch$-algebras equipped with a unit 1 whose compatibility with operations are defined as follows:
$$ x \succ 1 = x = 1 \succ x, \ \ \ x*1=x=1*x. $$
For instance, $Hoch_+(V):= K.1_K \oplus Hoch(V)$, where $Hoch(V)$ is the free $Hoch$-algebra over a $K$-vector space $V$ is a unital $Hoch$-algebra with unit $1_K$. This gives birth to unital $Mag^{\infty}$-algebras which are $Mag^{\infty}$-algebras such that the generating operations are related with the unit as follows:
$$ [1,\cdot, \ldots,  \cdot]_n=0,$$
$$ [\cdot, \ldots, 1, \ldots, \cdot]_n=[\cdot, \ldots, \cdot]_{n-1},$$
$$[\cdot,  \ldots, \cdot,1]_n=0.$$
Over $Hoch_+(V)$, one has a unital infinitesimal coproduct $\delta$ defined via the former coproduct $\Delta$ as follows:

$$ \delta(x) =1_K \otimes x+ x \otimes 1_K +\Delta(x), $$
for any $x \in Hoch(V)$.
The compatibility relations are the so-called
unital infinitesimal relations defined as follows:
$$ \Delta(x \succ y):= x_{(1)}\otimes (x_{(2)} \succ y) + (x \succ y_{(1)}) \otimes y_{(2)} - x\otimes y.$$
$$ \Delta(x * y):= x_{(1)}\otimes (x_{(2)} * y) + (x * y_{(1)}) \otimes y_{(2)} - x\otimes y.$$
We then obtain the good triple of operads $(As,Hoch,Mag^{\infty})$ equipped with the unital infinitesimal relations.
\section{Other triples of operads}
The triple of operads $(As,Hoch, Mag^\infty)$ endowed with the infinitesimal relations are not the only one involving the operad $Hoch$. By changing the compatibility relations, two other good triples of operads $(Com, Hoch, Prim_{Com}\ Hoch)$ and $(As, Hoch, Prim_{As}\ Hoch)$ endowed respectively with the Hopf relations and the semi-Hopf relations can be proposed. But contrary to the case of the triple $(As,Hoch, Mag^\infty)$
the explicit descriptions of operads of the primitive elements of these two other triples are open problems.

\bibliographystyle{plain}
\bibliography{These}

\end{document}